\documentclass[english]{article}
\usepackage[T1]{fontenc}
\usepackage{a4}
\usepackage{geometry}
\usepackage{amsmath}
\usepackage{babel}
\usepackage{amssymb}

\makeatletter

\providecommand{\LyX}{L\kern-.1667em\lower.25em\hbox{Y}\kern-.125emX\@}

 \newcommand{\lyxaddress}[1]{
   \par {\raggedright #1 
   \vspace{1.4em}
   \noindent\par}
 }

\newcommand{\meas}{\operatorname{meas}}

\makeatother
\begin{document}

\title{On a nonlinear compactness lemma in \( L^{p}(0,T;B) \).}

\author{\textbf{Emmanuel Maitre}\\
\\
\emph{Laboratoire de Mathématiques et Applications}\\
\emph{Université de Haute-Alsace}\\
\emph{4, rue des Frères Lumière}\\
\emph{68093 Mulhouse}\\
\texttt{\emph{\normalsize E.Maitre@uha.fr}}}

\maketitle
\begin{abstract}
We consider a nonlinear counterpart of a compactness lemma of J. Simon
\cite{Simon}, which arises naturally in the study of doubly nonlinear
equations of elliptic-parabolic type. Our work was motivated by previous
results J. Simon \cite{Simon}, recently sharpened by H. Amann \cite{Amann},
in the linear setting, and by a nonlinear compactness argument of
H.W. Alt and S. Luckhaus \cite{Alt}. \\
\\
\textbf{MSC2000 :} \emph{Primary 46B50, 47H30. Secondary 34G20, 35K65.}
\end{abstract}

\section{Introduction}

\textsl{\emph{Typical applications where the compactness argument
stated below is useful are those in which the following kind of doubly
nonlinear equations arises\[
\frac{dB(u)}{dt}+A(u)=f\]
 where \( A \) is elliptic and \( B \) monotone (not strictly).
It is the case, for example, in porous medium, semi-conductor equations,
... }}

\textsl{\emph{In our application, we considered the injection moulding
of a thermoplastic, with a mold of small thickness with respect to
its other dimensions. By averaging Navier-Stokes equations across
the thickness of the mold, and under an assumption (of Hele-Shaw)
stating that the velocity field is proportional to the pressure gradient,
the pressure equation can be written as a doubly nonlinear equation
\cite{6}. }}

\textsl{\emph{Note that in this context, the equation can degenerate
to an elliptic one. In order to get existence of a solution, one usually
perform a time discretization, use some result on elliptic operator
and pass to the limit as the time step goes to zero. In nonlinear
problems compactness in time and space is then required. The compactness
in space is easily obtained for \( u \) from a coerciveness assumption
on the elliptic part \( A \), but we have no estimate on \( \frac{\partial u}{\partial t} \)
since \( B \) could degenerate. Theorem 1 uses the space compactness
of \( u \) and some time regularity on \( B(u) \) to derive a compactness
for \( B(u) \), which in turn can be useful to pass to the limit
in nonlinear terms of \( A \) (provided \( A \) has a an appropriate
structure, e.g. \( B- \)pseudomonotone \cite{7}).}}

\section{Main result}

Let us consider two Banach spaces \( E_{1} \), \( E_{2} \). Let
\( T>0 \), \( p\in [1,+\infty ] \), and \( B \) a (nonlinear) \emph{compact}
operator from \( E_{1} \) to \( E_{2} \), i.e. which maps bounded
subsets of \( E_{1} \) to relatively compact subsets of \( E_{2} \).

\paragraph*{Theorem 1 :}

\emph{Let \( U \) be a bounded subset of \( L^{1}(0,T;E_{1}) \)
such that \( V=B(U) \) is a subset of \( L^{p}(0,T;E_{2}) \) bounded
in \( L^{r}(0,T;E_{2}) \) with \( r>1 \).} \textit{Assume}\textit{\emph{\begin{equation}
\label{time}
\lim _{h\to 0^{+}}\left\Vert v(\cdot +h)-v\right\Vert _{L^{p}(0,T-h;E_{2})}=0\qquad \textit {uniformly\, for\, }v\in V.
\end{equation}
}} \emph{Then \( V \) is relatively compact in \( L^{p}(0,T;E_{2}) \)}
(\emph{and in} \( \mathcal{C}(0,T;E_{2}) \) \emph{if} \( p=+\infty  \)).

\paragraph*{Remarks :}

\begin{enumerate}
\item \textsl{One can easily check that theorem 1 holds if we assume only}
\textsl{\emph{\( U \)}} \textsl{bounded in} \textsl{\emph{\( L^{1}_{loc}(0,T;E_{1}) \)}}
\textsl{and} \textsl{\emph{\( V \)}} \textsl{bounded in} \textsl{\emph{\( L^{r}_{loc}(0,T;E_{2}) \).}}
\item \textsl{In the case where \( B \) is the canonical injection from
\( E_{1} \) to \( E_{2} \), the assumption on \( B \) corresponds
to the compactness of the embedding of \( E_{1} \) into \( E_{2} \),
and the conclusion falls in the scope of previous results of J. Simon
\cite{Simon}, theorem 3.}
\item \textsl{The point in this result is that we do not make any structural
assumption on \( B \) (e.g. strict monotony, which would fall in
the scope of results of A. Visintin \cite{Visintin}) except compactness.
Note that in the case of a compact embedding of \( E_{1} \) into
\( E_{2} \), \( B \) needs only to be continuous from \( E_{1} \)
to \( E_{2} \) for the \( E_{2} \) topology.}
\end{enumerate}

\paragraph*{Idea of the proof :}

A sufficient condition for compactness is to prove that for each couple
\( (t_{1,}t_{2}), \) \( \int _{t_{1}}^{t_{2}}v(t)dt \) describes
a relatively compact subset of \( E_{2} \) as \( v \) describes
\( V \). First the \( u(t) \), \( u\in U \) are truncated in norm
at height \( M>0 \) and form a bounded subset of \( E_{1} \) which
\( B \) maps to a relatively compact subset \( V^{M}(t) \) of \( E_{2} \).
The key point is that thanks to equi-integrability assumption, \( \int _{t_{1}}^{t_{2}}v(t)dt \)
can be approximated uniformly in \( v \) by Riemann sums \emph{involving
truncated elements} of the \( V^{M}(t) \).

\paragraph*{Proof :}

Thanks to the equi-integrability (\ref{time}) of \( V \) and results
of \cite{Simon}, we only have to prove that for each \( (t_{1},t_{2}) \)
such that \( 0<t_{1}<t_{2}<T \), the set\[
K=\left\{ \int _{t_{1}}^{t_{2}}v(t)dt,\quad v\in V\right\} \]
is relatively compact in \( E_{2} \). For that purpose, we introduce
for \( u\in U \) and \( M>0 \) the measurable subset of \( [0,T] \)
defined by\[
G_{u}^{M}=\left\{ t\in [0,T],\quad \left\Vert u(t)\right\Vert _{E_{1}}>M\right\} .\]
From our assumptions on \( U \), there exists a constant \( C>0 \)
such that \[
\forall u\in U,\qquad \left\Vert u\right\Vert _{L^{1}(0,T;E_{1})}\leq C,\]
 and since we have \[
\operatorname {meas}(G_{u}^{M})=\int _{G_{u}^{M}}1dt\leq \int _{G^{M}_{u}}\frac{\left\Vert u(t)\right\Vert _{E_{1}}}{M}dt\leq \frac{C}{M}\]
 that gives\begin{equation}
\label{G->0}
\lim _{M\to +\infty }\operatorname {meas}(G_{u}^{M})=0,\qquad \textrm{ uniformly in }u.
\end{equation}
 Introducing the truncated functions\[
u^{M}(t)=u(t)\textrm{ if }t\not \in G^{M}_{u},\quad 0\textrm{ otherwise},\]
 we have by construction \begin{equation}
\label{bound}
\forall M>0,\quad \forall u\in U,\quad \forall t\in [0,T],\qquad \left\Vert u^{M}(t)\right\Vert _{E_{1}}\leq M.
\end{equation}

\paragraph*{Lemma 1}

\emph{Under condition (\ref{time}), \( K \) can be uniformly approximated
by Riemann sums involving elements of the form \( v^{M}(t)=B(u^{M}(t)) \),
in the following sense : given \( \varepsilon >0 \), there exist
integers \( N \) and \( M \) such that for all \( v=B(u)\in V \),
there exists \( s_{v}^{N,M}\in ]0,h[ \) such that \begin{equation}
\label{riem}
\left\Vert \int _{t_{1}}^{t_{2}}v(t)dt-\sum _{i=1}^{N}hv^{M}(\xi _{i-1}^{N}+s_{v}^{N,M})\right\Vert _{E_{2}}<\varepsilon 
\end{equation}
 where \( h=\frac{t_{2}-t_{1}}{N} \) and \( \xi _{i}^{N}=t_{1}+ih \). }

\paragraph*{Proof :}

{\raggedright We first note that \par}

\begin{equation}
\label{bid}
\int _{t_{1}}^{t_{2}}v(t)dt-\sum _{i=1}^{N}hv^{M}(\xi _{i-1}^{N}+s_{v}^{N,M})=\int _{t_{1}}^{t_{2}}\left( v(t)-\sum _{i=1}^{N}v^{M}(\xi _{i-1}^{N}+s_{v}^{N,M})\chi _{]\xi _{i-1}^{N},\xi _{i}^{N}]}(t)\right) dt.
\end{equation}
 Then we prove the following inequality, where \( r' \) stands for
the conjuguate exponent of \( r \) :

\begin{multline}\label{truc}\frac{1}{h}\int _{0}^{h}\int _{t_{1}}^{t_{2}}\left\Vert v(t)-\sum _{i=1}^{N}v^{M}(\xi _{i-1}^{N}+s)\chi _{]\xi _{i-1}^{N},\xi _{i}^{N}]}(t)\right\Vert _{E_{2}}dtds \\ \leq 2T^{1-\frac{1}{p}}\sup _{\sigma \in [-h,h]}\left\Vert v(\cdot +\sigma )-v\right\Vert _{L^{p}(0,T-\sigma ;E_{2})}+2\left(\meas G^M_u\right)^\frac{1}{r'}\left\|v-B(0)\right\|_{L^r(0,T;E_2)}.\end{multline}

{\raggedright Let us denote by \( I \) the left-hand side of the
stated inequality. Then\[
I=\frac{1}{h}\int _{0}^{h}\sum _{i=1}^{N}\int _{\xi _{i-1}^{N}}^{\xi _{i}^{N}}\left\Vert v(t)-v^{M}(\xi _{i-1}^{N}+s)\right\Vert _{E_{2}}dtds=\frac{1}{h}\sum _{i=1}^{N}\int _{\xi _{i-1}^{N}}^{\xi _{i}^{N}}\int _{\xi _{i-1}^{N}}^{\xi _{i}^{N}}\left\Vert v(t)-v^{M}(s)\right\Vert _{E_{2}}dtds.\]
Using Fubini's theorem, and setting \( \sigma =s-t \) we get \[
I=\frac{1}{h}\sum _{i=1}^{N}\int _{\xi _{i-1}^{N}}^{\xi _{i}^{N}}\int _{\xi _{i-1}^{N}-t}^{\xi _{i}^{N}-t}\left\Vert v(t)-v^{M}(t+\sigma )\right\Vert _{E_{2}}d\sigma dt,\]
 which gives thanks to a new application of Fubini's theorem,\[
I=\frac{1}{h}\int _{-h}^{h}\sum _{i=1}^{N}\int _{max(\xi _{i-1}^{N},\xi _{i-1}^{N}-\sigma )}^{min(\xi _{i}^{N},\xi _{i}^{N}-\sigma )}\left\Vert v(t)-v^{M}(t+\sigma )\right\Vert _{E_{2}}dtd\sigma \leq \frac{1}{h}\int _{-h}^{h}\int _{max(t_{1},t_{1}-\sigma )}^{min(t_{2},t_{2}-\sigma )}\left\Vert v(t)-v^{M}(t+\sigma )\right\Vert _{E_{2}}dtd\sigma .\]
From the definition of \( v^{M} \) we thus have\[
I\leq \frac{1}{h}\int _{-h}^{h}\int _{max(t_{1},t_{1}-\sigma )}^{min(t_{2},t_{2}-\sigma )}\left\Vert v(t)-v(t+\sigma )\right\Vert _{E_{2}}dtd\sigma +\frac{1}{h}\int _{-h}^{h}\int _{max(t_{1},t_{1}-\sigma )}^{min(t_{2},t_{2}-\sigma )}\chi _{G_{u}^{M}}(t+\sigma )\left\Vert v(t)-B(0)\right\Vert _{E_{2}}dtd\sigma .\]
 As \( V \) is a bounded subset of \( L^{r}(0,T;E_{2}) \) one has
the second term bounded by\par}

\[
\frac{1}{h}\int _{-h}^{h}\left( \int _{\max (t_{1},t_{1}-\sigma )}^{\min (t_{2},t_{2}-\sigma )}\chi _{G_{u}^{M}}(t+\sigma )dt\right) ^{\frac{1}{r'}}\left( \int _{t_{1}}^{t_{2}}||v(t)-B(0)||_{E_{2}}^{r}dt\right) ^{\frac{1}{r}}d\sigma \leq 2(\meas G_{u}^{M})^{\frac{1}{r'}}||v-B(0)||_{L^{r}(0,T;E_{2})}.\]

{\raggedright and the Hölder inequality gives the announced estimation
(\ref{truc}).\par}

Using (\ref{time}), (\ref{G->0}) and as \( v \) belongs to a bounded
subset \( V \) of \( L^{r}(0,T;E_{2}) \), we conclude from (\ref{truc})
that\begin{equation}
\label{unif1}
\frac{1}{h}\int _{0}^{h}\int _{t_{1}}^{t_{2}}\left\Vert v(t)-\sum _{i=1}^{N}v^{M}(\xi _{i-1}^{N}+s)\chi _{]\xi _{i-1}^{N},\xi _{i}^{N}]}(t)\right\Vert _{E_{2}}dtds\to 0,\qquad \textrm{when }M\textrm{ and }N\textrm{ go to infinity},\textrm{ uniformly in }v.
\end{equation}
 We claim that there exists at least one \( s=s_{v}^{N,M}\in [0,h] \)
such that \begin{equation}
\label{unif2}
\int _{t_{1}}^{t_{2}}\left\Vert v(t)-\sum _{i=1}^{N}v^{M}(\xi _{i-1}^{N}+s_{v}^{N,M})\chi _{]\xi _{i-1}^{N},\xi _{i}^{N}]}(t)\right\Vert _{E_{2}}dt\to 0,
\end{equation}
when \( M,N \) go to infinity, uniformly in \( v. \) Indeed, let
us set by sake of readability\[
f^{v}_{N,M}(s)=\int _{t_{1}}^{t_{2}}\left\Vert v(t)-\sum _{i=1}^{N}v^{M}(\xi _{i-1}^{N}+s)\chi _{]\xi _{i-1}^{N},\xi _{i}^{N}]}(t)\right\Vert _{E_{2}}dt\]
 so that the uniform convergence (\ref{unif1}) reads \begin{equation}
\label{bidule}
\frac{1}{h}\int _{0}^{h}f^{v}_{N,M}(s)ds\to 0,\qquad \textrm{when }M\textrm{ and }N=\frac{1}{h}\textrm{ go to infinity},\textrm{ uniformly in }v.
\end{equation}
Then for fixed \( v,N,M \) there exists at least one \( s=s_{v}^{N,M}\in [0,h] \)
such that\[
f^{v}_{N,M}(s_{v}^{N,M})\leq \frac{1}{h}\int _{0}^{h}f^{v}_{N,M}(s)ds.\]
If not, we would have the reverse strict inequality for all \( s\in [0,h] \)
which by averaging on \( [0,h] \) would lead to a contradiction.
Then as \( f_{N,M}^{v} \) is positive, the uniform convergence (\ref{bidule})
implies 

\begin{equation}
f_{N,M}^{v}(s_{v}^{N,M})\to 0,\qquad \textrm{when }M\textrm{ and }N=\frac{1}{h}\textrm{ go to infinity},\textrm{ uniformly in }v,
\end{equation}
 which is exactly (\ref{unif2}).

\emph{A fortiori}, (\ref{riem}) holds thanks to (\ref{bid}) and
since\[
\left\Vert \int _{t_{1}}^{t_{2}}\left( v(t)-\sum _{i=1}^{N}v^{M}(\xi _{i-1}^{N}+s_{v}^{N,M})\chi _{]\xi _{i-1}^{N},\xi _{i}^{N}]}(t)\right) dt\right\Vert _{E_{2}}\leq \int _{t_{1}}^{t_{2}}\left\Vert v(t)-\sum _{i=1}^{N}v^{M}(\xi _{i-1}^{N}+s_{v}^{N,M})\chi _{]\xi _{i-1}^{N},\xi _{i}^{N}]}(t)\right\Vert _{E_{2}}dt\]
 This proves lemma 1. To conclude the proof of theorem 1, note that
lemma 1 means that \( K\subset \varepsilon B_{E_{2}}+K_{M,N} \) where
\( B_{E_{2}} \) is the unit open ball of \( E_{2} \) and \[
K_{M,N}=\left\{ \sum _{i=1}^{N}hv^{M}(\xi _{i-1}^{N}+s_{v}^{N,M}),\quad v^{M}=B(u^{M}),\quad u\in U\right\} .\]
For fixed \( M,N \) and from (\ref{bound}) we note that \( u^{M}(\xi _{i-1}^{N}+s_{v}^{N,M}) \)
is bounded in \( E_{1} \) uniformly in \( u\in U \). As \( B \)
is compact, \( K_{M,N} \) is thus a relatively compact subset of
\( E_{2} \). Thus \( K \) is also relatively compact in \( E_{2} \).
\( \diamondsuit  \)

\paragraph*{Corollary 1 :}

\emph{Let \( U \) be a bounded subset of \( L^{1}(0,T;E_{1}) \)
such that \( V=B(U) \) is bounded in \( L^{r}(0,T;E_{2}) \) with
\( r>1 \).} \textit{Assume\[
\frac{\partial V}{\partial t}=\left\{ \frac{\partial v}{\partial t},\quad v\in V\right\} \textrm{ is bounded in }L^{1}(0,T;E_{2}).\]
} \emph{Then \( V \) is relatively compact in \( L^{p}(0,T;E_{2}) \)
for any} \( p<+\infty  \).

\paragraph*{Proof : }

Condition (\ref{time}) of theorem 1 is satisfied (see \cite{Simon},
lemma 4).

\bigskip{}
\lyxaddress{\emph{The author would like to thanks the referee for its careful
reading of his paper, which lead to an improved formulation of the
results.}}

\end{document}